\documentclass[preprint]{elsarticle}       % Specifies the document class

\usepackage{array}

\usepackage{lineno,hyperref}
\usepackage{rotating}
\usepackage{makeidx}
\usepackage{float}
\usepackage{epsfig}
\usepackage{amsmath,amssymb,amsthm}
\usepackage[latin1]{inputenc}
\usepackage{epic,eepic}
\usepackage{overpic}
\usepackage{color}
\usepackage{amsfonts}
\usepackage{graphicx}
\usepackage{mwe} % For dummy images
\usepackage{subcaption}

\usepackage{algorithm}
\usepackage{algorithmicx}
\usepackage{algpseudocode}

            {\hfill\penalty10000\raisebox{-.09em}{\large\bf\rm $\blacksquare$}\par\medskip}

\newtheorem{theorem}{Theorem}[section]

\newtheorem{corollary}[theorem]{Corollary}

\newtheorem{remark}{Remark}[section]
\begin{document}

\small

\begin{frontmatter}

\title{{\bf A general correction for numerical integration rules over piece-wise continuous functions.}}

\author[IIPE]{Shipra Mahata}
\ead{shipramahata.maths@iipe.ac.in}
\author[IIPE]{Samala Rathan}
\ead{rathans.math@iipe.ac.in}
\author[UPCT]{Juan Ruiz-\'Alvarez}
\ead{juan.ruiz@upct.es}
\author[UV]{Dionisio F. Y\'an\~ez}
\ead{dionisio.yanez@uv.es}
\date{Received: date / Accepted: date}

\address[IIPE]{Department Mathematics, Indian Institute of Petroleum and Energy (IIPE) Visakhapatnam, Andhra Pradesh, India.}
\address[UPCT]{Departamento de Matem\'atica Aplicada y Estad\'istica. Universidad  Polit\'ecnica de Cartagena. Cartagena, Spain.}
\address[UV]{Departamento de Matem\'aticas. Universidad de Valencia. Valencia, Spain.}

%The correct dates will be entered by the editor.

\begin{abstract}
This article presents a novel approach to enhance the accuracy of classical quadrature rules by incorporating correction terms. The proposed method is particularly effective when the position of an isolated discontinuity in the function and the jump in the function and its derivatives at that position are known. Traditional numerical integration rules are exact for polynomials of certain degree. However, they may not provide accurate results for piece-wise polynomials or functions with discontinuities without modifying the location and number of data points in the formula. Our proposed correction terms address this limitation, enabling the integration rule to conserve its accuracy even in the presence of a jump discontinuity. The numerical experiments that we present support the theoretical results obtained.%This advancement in numerical integration techniques opens up new possibilities in the study of this field.% for more accurate and efficient computation of in a variety of scientific and engineering applications.
\end{abstract}
\end{frontmatter}
%All acknowledgements should be placed in the back of the paper after Conclusions..

\section{Introduction}
Numerical integration is a fundamental tool in many areas of science and engineering \cite{book1, book2, book3, book4, Gautschi, Tref, Tref1}. Classical quadrature rules, such as Gaussian quadrature, provide efficient and accurate methods for approximating definite integrals. However, these methods can lose accuracy when applied to piece-wise polynomials or functions with finite discontinuities \cite{Amat2023}. 

%Classical integration formulas, such as the trapezoidal rule, the Simpson's rule, or the Newton-Cotes formulas, are based on the integration of interpolatory polynomials over an interval. The classical problem that arises from using such interpolatory polynomials is the loss of accuracy whenever the original data does not present enough regularity. In this article, we introduce a new method that consists in the construction and use of a polynomial correction that takes into account the presence of discontinuities in the integrand. The problem of obtaining quadrature rules adapted to the presence of discontinuities in this context can be found in the literature \cite{GRIER2014193, Tornberg2002}, but we have not found many references about the subject. 

Conventional integration formulas, such as the trapezoidal rule, Simpson's rule, and Newton-Cotes formulas are fundamentally based on the integration of interpolatory polynomials over a specified interval. A usual problem that arises from the use of such interpolatory polynomials is the degradation of accuracy when the original data do not present enough regularity. This problem also affects other numerical integration techniques, such as Gaussian quadrature formulas.

In this article, we propose a novel method that involves the construction and application of a polynomial correction that accounts for the presence of discontinuities in the function that needs to be integrated. This method aims to provide a solution for the problem of the presence of finite discontinuities in the integrand when trying to apply classical numerical integration formulas. The challenge of developing quadrature rules that are adapted to the presence of discontinuities in this context has been discussed in the literature \cite{GRIER2014193, Tornberg2002, libro_shizgal, Amat2023}. However, it is also true that we have not found many references that discuss about this topic. Our aim is to introduce a novel approach to address the mentioned issue by incorporating correction terms into the classical quadrature rules. The proposed correction terms have an explicit form and are easy to compute if the position of the discontinuity and the jump in the function and its derivatives at that position are known. This is a significant advantage over other methods, which often require complex computations or approximations. %The correction terms allow the integration rule to be exact for piece-wise polynomials of the degree, 
The correction terms allow the quadrature rule to reproduce piece-wise polynomials of the same degree that the integration formula reproduces at smooth areas, enhancing the accuracy of the numerical integration.
\par
This technique is of special interest for Gauss-Legendre quadrature formulas that use a large number of points. In the presence of a discontinuity, the typical adaptation would be to split the interval of integration in two, effectively doubling the number of points. However, this can significantly increase the computational cost. By using the proposed correction terms, we can maintain the accuracy of the integration without needing to increase the number of points.
\par
The proposed method is general and enhances the accuracy of classical quadrature rules in the presence of discontinuities, providing an efficient and easy-to-implement solution for the numerical integration of piece-wise smooth functions. %We believe this advancement opens up new possibilities for more accurate and efficient computation in a variety of scientific and engineering applications.
\par
This article is organized as follows: Section \ref{corr_terms_general} presents how to obtain a general polynomial correction term for functions with an isolated discontinuity in the interval of integration without modifying the number of points used in the quadrature formula. We also analyze and provide a formula for the error involved in the numerical integration when the location of a discontinuity is known. In Section \ref{par_gauss}, we discuss the particularization of the proposed polynomial correction to the case of Gaussian quadrature. In Section \ref{acc_x*} we discuss the accuracy needed in the location of the discontinuity and in the jumps in the function and its derivatives if these quantities are to be approximated numerically. Section \ref{num_exp} presents an extended battery of numerical experiments. Finally, Section \ref{conc} presents some conclusions.

\section{The correction term for a general case}\label{corr_terms_general}

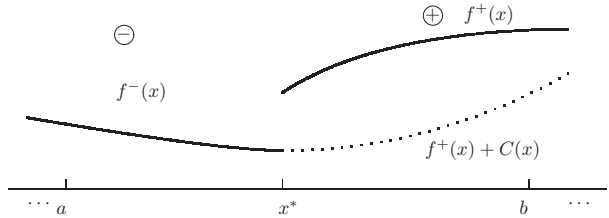
\begin{figure}[ht]
\begin{center}
\resizebox{13cm}{!} {
\begin{picture}(500,110)(0,0)
%escala k=2
\put(50,-10){\line(1,0){310}}
%\put(0,-10){\line(0,1){5}}
\put(80,-10){\line(0,1){5}}
%\put(160,-10){\line(0,1){5}}
%\put(200,-10){\line(0,1){5}}%x_{j+1/2}
\put(192,-10){\line(0,1){5}}%x^*
%\put(240,-10){\line(0,1){5}}
\put(320,-10){\line(0,1){5}}
%\put(400,-10){\line(0,1){5}}
%\put(480,-10){\line(0,1){5}}
%celdas escala k=2
\put(60,-20){$\cdots$}
%\put(0,86){$f_{j-2}$}
%\put(80,30){$f_{j-1}$}
%\put(160,23){$f_{j}$}
%\put(235,74){$f_{j+1}$}
%\put(315,83){$f_{j+2}$}
%\put(400,186){$f_{j+3}$}
\put(340,-20){$\cdots$}

%\put(-5,-23){$x_{j-2}$}
\put(75,-23){$a$}
%\put(155,-23){$x_{j}$}
%\put(195,-23){$x_{j+1/2}$}
\put(190,-23){$x^*$}
%\put(235,-23){$x_{j+1}$}
\put(315,-23){$b$}
%\put(395,-23){$x_{j+3}$}
%\put(475,-23){$x_{j+4}$}

%86    30     6    74   114   

%\put(0,-10){\line(0,1){86}}
%\put(0,76){\circle{5}}
%\put(80,-10){\line(0,1){33}}
%\put(80,23){\circle{5}}
%\put(160,-10){\line(0,1){24}}
%\put(160,14){\circle{5}}
%\put(240,-10){\line(0,1){72}}
%\put(240,62){\circle{5}}
%\put(320,-10){\line(0,1){83}}
%\put(320,73){\circle{5}}
%\put(400,-10){\line(0,1){186}}
%\put(400,176){\circle{5}}

%polinomio izquierda
\linethickness{0.3mm}
\qbezier(60,27)(160,10)(192,10)
%polinomio derecha
\qbezier(192,40)(240,72)(340,73)
%polinomio completo punteado
\qbezier[30](192,10)(270,10)(340,50)

%lado menos
\put(110,70){\circle{10}}
\put(106,68){$-$}
\put(106,38){$f^-(x)$}
%lado mas
\put(270,80){\circle{10}}
\put(266,78){$+$}
\put(286,78){$f^+(x)$}
\put(266,8){$f^+(x)+C(x)$}
\end{picture}
}
\end{center}
\caption{In this figure we can see an example of an isolated discontinuity in the function and its derivatives placed in the interval $[a,b]$ at a position $x^*$.}\label{disc_alpha_corner}
\end{figure}

As shown in Figure \ref{disc_alpha_corner}, our initial assumption is that the location $x^*$ of an isolated discontinuity, which lies within the interval $[a,b]$, is either known or can be estimated with a high degree of accuracy (in Section \ref{acc_x*} we discuss the accuracy needed in the location of $x^*$). Then, we can explain how to incorporate a correction term into a general quadrature formula to allow for the presence of the discontinuity. 

%===
Let $l$ be a natural number and $f^-,f^+\in\mathcal{C}^{l+1}([a,b])$ be two functions, where the superindex $-$ denotes the function values to the left of the discontinuity, and the superindex $+$ denotes the function values to the right of the discontinuity. The function $f(x)$ is then defined as,
\begin{equation}\label{eq1}
f(x)=\begin{cases}
f^-(x),& \text{for } a\le x< x^*,\\
f^+(x),& \text{for } x^*\leq x \leq b.
\end{cases}
\end{equation}
We express the discontinuities in the derivatives of $f$ as,
\begin{equation}\label{ir}
\left[f^{(k)}\right]=\lim_{x\to x^*+}f_{\underbrace{x\cdots x}_k}^{+}(x) - \lim_{x\to x^*-}f_{\underbrace{x\cdots x}_k}^{-}(x),
\end{equation}
where $f^{\pm}_{\underbrace{x\cdots x}_k}(x)$ represents the $k^{th}$ derivative of $f^\pm$ at $x$. %The approximation of these discontinuities and the required accuracy will be discussed later.

%===
%
%
%Let us suppose that we want to integrate the piece-wise continuous function in the interval $[-1,1]$,
%\ 
%begin{equation}\label{eq1}
%%\scalemath{0.9}{
%f(x)=\left\{\begin{array}{ll}
%f^-(x),  & \textrm{ if } -1\le x<x^*,\\
%f^+(x),  & \textrm{ if } x^*\le x\le 1,
%\end{array}
%\right.
%%}
%\end{equation}
\begin{theorem}\label{teo1}
Given the piece-wise continuous function $f(x)$ described previously in \eqref{eq1}, we can write its integral in the interval $[a,b]$ as,
\begin{equation}\label{int1kk}
\begin{aligned}%\nonumber
\int_{a}^{b}f(x)dx&=\int_{a}^{b}f^-(x)dx+\int_{x^*}^{b}C(x)dx+E,
\end{aligned}
\end{equation}
where,
\begin{equation}\label{intC1kk}
\begin{aligned}%\nonumber
\int_{x^*}^{b}C(x)dx&=\sum_{k=1}^{l+1}\frac{\left[f^{(k)}\right]}{k!}(b-x^*)^k,
\end{aligned}
\end{equation}
being the truncation error of the integral $E$ bounded by
%\begin{equation}\label{errteo}
%E=\int_{x^*}^b \frac{(f^+)^{(l+1)}(\xi_1(x))-(f^-)^{(l+1)}(\xi_2(x))}%{(l+1)!}(x-x^*)^{l+1},\quad \xi_i(x)\in(x^*,x), \, \, i=1,2.
%\end{equation}
\begin{equation}\label{errteo}
E\leq \frac{||(f^+)^{(l+1)}||_{\infty,[a,b]}+||(f^-)^{(l+1)}||_{\infty,[a,b]}}{(l+1)!}(b-x^*)^{l+2}.
\end{equation}

\end{theorem}
\begin{proof}
We want to integrate $f(x)$. Let us develop the expressions up to the index $l$ using Taylor's expansions, so that our formulas are exact for polynomials of degree smaller or equal than $l$. In this case, we can write, 
\begin{equation}\label{fmmi1}
f^+(x)-f^-(x)=\sum_{k=0}^{l}\frac{[f^{(k)}]}{k!}(x-x^*)^k+e(x), \quad x\ge x^*,
%[f]+[f'](x-x^*)+\frac{[f'']}{2!}(x-x^*)^2+\frac{[f''']}{3!}(x-x^*)^3+\cdots+\frac{[f^{(n)}]}{n!}(x-x^*)^n+E(x), \quad x\ge x^*,
\end{equation}
where we make use of the jump conditions in (\ref{ir}). Of course, $e(x)$ represents the local truncation error of the polynomial approximation in (\ref{fmmi1}), that can be written as,
\begin{equation}\label{ex}
e(x)=\frac{(f^+)^{(l+1)}(\xi_1(x))-(f^-)^{(l+1)}(\xi_2(x))}{(l+1)!}(x-x^*)^{l+1},\quad \xi_i(x)\in(x^*,x), \, \, i=1,2,
\end{equation}
and that is equal to zero if $f(x)$ is a piece-wise polynomial of degree smaller or equal than $l$.
%\begin{equation}\label{ir}
%\begin{aligned}%\nonumber
%\left[f\right]&=\lim_{x\to x^*}\left(f^{+}(x) - f^{-}(x)\right),\\
%\left[ f' \right]&=\lim_{x\to x^*}\left(f_x^{+}(x) - f_x^{-}(x)\right),\\
%\left[ f'' \right]&=\lim_{x\to x^*}\left(f_{xx}^{+}(x) - f_{xx}^{-}(x)\right),\\
%\left[ f''' \right]&=\lim_{x\to x^*}\left(f_{xxx}^{+}(x) - f_{xxx}^{-}(x)\right).
%\end{aligned}
%\end{equation}
%In other case, $E(x)\neq 0$ {\color{red}: can we prove anything else about the error?}.
Let us denote 
\begin{equation}\label{C}
C(x)=\sum_{i=0}^{l}\frac{[f^{(i)}]}{i!}(x-x^*)^i, \quad x\ge x^*,%+E(x)
%[f]+[f'](x-x^*)+\frac{[f'']}{2!}(x-x^*)^2+\frac{[f''']}{3!}(x-x^*)^3+\cdots+\frac{[f^{(n)}]}{n!}(x-x^*)^n+E(x), \quad x\ge x^*,
\end{equation}
so that from (\ref{fmmi1}) we can write,
\begin{equation}\label{fmmi2}
f^+(x)=f^-(x)+C(x)+e(x),\quad x\ge x^*.
\end{equation}
Taking into account these considerations, it should be clear that we can write the definite integral of $f(x)$ in $[a,b]$ as,

\begin{equation}\label{int1}
\begin{aligned}%\nonumber
\int_{a}^{b}f(x)dx&=\int_{a}^{x^*}f^-(x)dx+\int_{x^*}^{b}f^+(x)dx=\int_{a}^{x^*}f^-(x)dx+\int_{x^*}^{b}\left(f^-(x)+C(x)+e(x)\right)dx\\
&=\int_{a}^{b}f^-(x)dx+\int_{x^*}^{b}C(x)dx+\int_{x^*}^{b}e(x)dx.
\end{aligned}
\end{equation}
%The first part of the integral can be obtained using any classical quadrature formulas. 
The second part of the integral can be obtained explicitly using induction as,
\begin{equation}\label{intC1}
\begin{aligned}%\nonumber
\int_{x^*}^{b}C(x)dx=\sum_{i=1}^{l+1}\frac{[f^{(i)}]}{i!}(b-x^*)^i,
%[f] (b-x^*)+\frac{1}{2} [f'] (b-x^*)^2+\frac{1}{6} [f''] (b-x^*)^3+\frac{1}{24} [f'''] (b-x^*)^4+\int_{x^*}^bE(x)dx.
\end{aligned}
\end{equation}
and
$$E=\int_{x^*}^be(x)dx,$$
with $e(x)$ given in (\ref{ex}).
%\textcolor{blue}{Is there any general expression to this term?}
\end{proof}

\begin{corollary}
The truncation error $E$ is zero if $f(x)$ is a piece-wise polynomial of degree smaller or equal than $l$. 
\end{corollary}
\begin{proof}
The proof is direct just looking at the expression of the error (\ref{errteo}) presented in Theorem \ref{teo1}.
\end{proof}

\begin{remark}
The integral $\int_{a}^{b}f^-(x)dx$ in Theorem (\ref{teo1}) is susceptible of being computed using any quadrature rule. Thus, the number of terms included in (\ref{intC1}), i.e. $l$, should be adjusted so that the order of accuracy of the quadrature rule chosen for the integral  $\int_{a}^{b}f^-(x)dx$ is conserved. If we choose the Gaussian quadrature for the mentioned integral, then we should also analyse what happens when the weight function $w(x)$ is different from one. 
We think that the case of Gaussian quadrature is of special interest, as the expression of $f(x)$ must be known, including the position of the discontinuity. In the case of a piece-wise smooth function with an isolated discontinuity in the interval of integration, the easiest way of integrating would be to split the interval of integration in two. But in this case, the number of points used for obtaining an accurate result is multiplied by two, doubling the computational cost. This can be significant if the number of points used is high, or if the calculation involves several integrals of piece-wise smooth functions. We analyse the case of Gaussian quadrature in the next section.
\end{remark}

 \par
Let us try to get more insight about the use of the correction term proposed before in integrals of the form,  
\begin{equation}\label{eq2}
\int_{a}^{b}\omega(x)f(x)dx,
\end{equation}
with $\omega:[a,b]\to \mathbb{R}$ a continuous and positive function. In this case, if we suppose that the location of a discontinuity $x^*$ is known or approximated, we can propose different strategies to calculate this integral. One option is to divide the interval $[a,b]$ into two sub-intervals $[a,x^*]\cup[x^*,b]$ and apply an integration rule, for example, a Gaussian quadrature. As mentioned before, this strategy doubles the number of points considered in the computation. Another strategy can be to define the function 
\begin{equation}\label{eq3}
\widetilde{f}(x)=\begin{cases}
f^-(x),& \text{for } a\le x< x^*,\\
f^+(x)-T_{x^*}^l(f^+-f^-)(x),& \text{for } x^*\leq x \leq b,
\end{cases}
\end{equation}
being $T_{x^*}^l(f^+-f^-)$ the Taylor polynomial of the function $(f^+-f^-)$ developed at point $x^*$, i.e. 
\begin{equation}\label{tayl_T}
T_{x^*}^l(f^+-f^-)(x)=\sum_{k=0}^{l}\frac{\left[f^{(k)}\right]}{k!}(x-x^*)^k.
\end{equation}
Note that $\widetilde{f}$ is well defined if $f^-,f^+\in\mathcal{C}^{l+1}([a,b])$. Also, $\widetilde{f}\in\mathcal{C}^{l}([a,b])$ since
$$\widetilde{f}^{(j)}(x^*)=(f^+)^{(j)}(x^*)-(T_{x^*}^l)^{(j)}(f^+-f^-)(x^*)=
(f^+)^{(j)}(x^*)-(f^+-f^-)^{(j)}(x^*)=(f^-)^{(j)}(x^*), \quad j=0,\hdots,l.
$$
%It is easy to check that if $x\in [x^*,b]$ then there exists %$\xi_x\in(x^*,x)$
%$$\widetilde{f}(x)=T_{x^*}^l(f^-)(x)+(f^+(x)-T_{x^*}^l(f^+)(x))=T_{x^*}^l(f^-)(x)+\frac{(f^+)^{l+1)}(\xi_x)}{(l+1)!}(x-x^*)^{l+1}.$$
Now we approximate \eqref{eq2} using an integration rule $\mathcal{S}^\omega$ over $\tilde{f}$ in the whole interval $[a,b]$ plus the approximation of the integral of the Taylor polynomial in $[x^*,b]$, that in some cases can be easy to calculate, for example if $\omega(x)=1$, or that can be approximated using an appropriate integration rule, that we denote by $\tilde{\mathcal{S}}^\omega_{[x^*,b]}(T^l_{x^*}(f^+-f^-))$
\begin{equation}\label{eq4}
\int_{a}^{b}\omega(x)f(x)dx\approx \mathcal{S}^\omega_{[a,b]}(\widetilde{f})+\int_{x^*}^b \omega(x)T_{x^*}^l(f^+-f^-)(x)dx\approx \mathcal{S}^\omega_{[a,b]}(\widetilde{f})+\tilde{\mathcal{S}}^\omega_{[x^*,b]}(T^l_{x^*}(f^+-f^-))=:\mathcal{N}^\omega_{[a,b]}(f).
\end{equation}
Note that, if $\omega(x)=1$, the correction term has a closed formula for all $x$.

If we analyze now the error for the proposed integration process, after some algebraic manipulations we get the same error obtained when we apply the integration rule over $\widetilde{f}$ plus another term,
\begin{equation}
\begin{split}
\int_{a}^{b}\omega(x)f(x)dx-\mathcal{N}^\omega_{[a,b]}(f)=&\int_{a}^{b}\omega(x)f(x)dx-\mathcal{S}^\omega_{[a,b]}(\widetilde{f})-\tilde{\mathcal{S}}^\omega_{[x^*,b]}(T^l_{x^*}(f^+-f^-))\\
=&\int_{a}^{b}\omega(x)f(x)dx-\mathcal{S}^\omega_{[a,b]}(\widetilde{f})-\int_{x^*}^b \omega(x)T_{x^*}^l(f^+-f^-)(x)dx\\
&+\int_{x^*}^b \omega(x)T_{x^*}^l(f^+-f^-)(x)dx-\tilde{\mathcal{S}}^\omega_{[x^*,b]}(T^l_{x^*}(f^+-f^-))\\
=&\int_{a}^{x^*}\omega(x)f(x)dx+\int_{x^*}^b\omega(x)f(x)dx-\mathcal{S}^\omega_{[a,b]}(\widetilde{f})-\int_{x^*}^bw(x) T_{x^*}^l(f^+-f^-)(x)dx\\
&+\int_{x^*}^b \omega(x)T_{x^*}^l(f^+-f^-)(x)dx-\tilde{\mathcal{S}}^\omega_{[x^*,b]}(T^l_{x^*}(f^+-f^-))\\
=&\int_{a}^{x^*}\omega(x)f^-(x)dx+\int_{x^*}^b\omega(x)\left(f^+(x)-T_{x^*}^l(f^+-f^-)(x)\right)dx-\mathcal{S}^\omega_{[a,b]}(\widetilde{f})\\
&+\int_{x^*}^b \omega(x)T_{x^*}^l(f^+-f^-)(x)dx-\tilde{\mathcal{S}}^\omega_{[x^*,b]}(T^l_{x^*}(f^+-f^-))\\
=&\int_{a}^{b}\omega(x)\widetilde{f}(x)dx-\mathcal{S}^\omega_{[a,b]}(\widetilde{f})+\int_{x^*}^b \omega(x)T_{x^*}^l(f^+-f^-)(x)dx-\tilde{\mathcal{S}}^\omega_{[x^*,b]}(T^l_{x^*}(f^+-f^-)).\\
\end{split}
\end{equation}
So we obtain that the error between the real integral of $f$ and our new approximation method can be expressed as the sum of two errors: one that comes from the approximation of the integral of $\tilde{f}$ and another one that comes from the approximation of the integral of the Taylor polynomial $T_{x^*}^l(f^+-f^-)(x)$. If this last error is equal to $0$ (because the integration rule is exact for polynomials of degree l), then the error of the proposed approximation is the same as the one obtained for the approximation of the integral of the continuous function $\tilde{f}$ in (\ref{eq3}).

\section{Particularisation to Gaussian quadrature}\label{par_gauss}

In the general case of Gaussian quadrature, we want to obtain formulas for the approximated evaluation of the integral
$$\int_a^b w(x)f(x)dx.$$
In the classical approach, $f$ is defined in $[a, b]$ and it is continuous and differentiable on this interval. Here we consider that $f$ might present an isolated discontinuity at $x^*$ in the interval $(a,b)$, that might affect the function or  the derivatives of $f$. On the other hand, $w$ is a weight function that is defined, positive, continuous and integrable on $(a,b).$ Taking into account the expressions (\ref{fmmi1}-\ref{fmmi2}),
%\begin{equation}\label{C}
%C(x)=\sum_{i=0}^{n}\frac{[f^{(i)}]}{i!}(x-x^*)^i+\sum_{i=n+1}^{\infty}\frac{[f^{(i)}]}{i!}(x-x^*)^i, \quad x\ge x^*,%+E(x)
%%[f]+[f'](x-x^*)+\frac{[f'']}{2!}(x-x^*)^2+\frac{[f''']}{3!}(x-x^*)^3+\cdots+\frac{[f^{(n)}]}{n!}(x-x^*)^n+E(x), \quad x\ge x^*,
%\end{equation}
%so that from (\ref{fmmi1}) we can write,
%\begin{equation}\label{fmmi2}
%f^+(x)=f^-(x)+C(x),\quad x\ge x^*.
%\end{equation}
we can write,
\begin{equation}\label{CG}
\begin{aligned}
\int_a^b w(x)f(x)dx=\int_a^{x^*} w(x)f^-(x)dx+\int_{x^*}^b w(x)\left(f^-(x)+C(x)+e(x)\right)dx=\int_a^b w(x)f^-(x)dx+\int_{x^*}^b w(x)C(x)dx+\int_{x^*}^b w(x)e(x)dx,
\end{aligned}
\end{equation}
that applying Gauss quadrature rule, can be written as,
\begin{equation}\label{CG2}
\begin{aligned}
\int_a^b w(x)f(x)dx\approx \sum_{k=0}^l W_k f^-(x_k)+\sum_{k=0}^l \tilde{W}_k C(x_k),
\end{aligned}
\end{equation}
where the quadrature weights have the expression,
\begin{equation}\label{w}
\begin{aligned}
W_k&=\int_a^b w(x)[L_k(x)]^2 dx,\\
\tilde{W}_k&=\int_{x^*}^b w(x)[\tilde{L}_k(x)]^2 dx,\\
L_k(x)&=\prod_{i=0;i\neq k}^l  \frac{x-x_i}{x_k-x_i}, k=0,1,\cdots, l,\\
\tilde{L}_k(x)&=\prod_{i=0;i\neq k}^l  \frac{x-\tilde{x}_i}{\tilde{x}_k-\tilde{x}_i}, k=0,1,\cdots, l,
\end{aligned}
\end{equation}
and $x_k,\tilde{x}_k, k=0,\ldots, l$, are the zeros of polynomials of degree $l+1$ from a system of orthogonal polynomials to $w(x)$ in $(a,b)$ or $(x^*, b)$ respectively. See Chapter 10 of \cite{book1} for more information about how the previous formulas are obtained.  Let us now analyse the error introduced by this approximation. We can adapt Theorem 10.1 of \cite{book1} to our case:

\begin{theorem}
Let us suppose that $w$ is a weight function that is defined, integrable, continuous and positive in the interval $(a, b)$, and that $f$ in (\ref{eq1}) is defined and piece-wise continuous on $[a,b]$ with an isolated discontinuity at $x^*$. Let us also suppose that $f^-$ has a continuous derivative of order $2l+2$ on $[a,b], l\ge0$. Then, it exists a number $\eta\in (a,b)$ that satisfies,
$$\int_a^b w(x)f(x)dx- \sum_{k=0}^l W_k f^-(x_k)-\sum_{k=0}^l \tilde{W}_k C(x_k)=K_l\cdot (f^-)^{(2l+2)}(\eta)+\int_{x^*}^bw(x)e(x) dx,$$
with
$$K_l=\frac{1}{(2l+2)!}\int_a^b w(x)[(x-x_0)\cdots(x-x_l)]^2dx.$$
So the integration formula (\ref{CG2}-\ref{w}) is exact for every piece-wise polynomial of degree $2l+1$ or less with an isolated discontinuity at $x^*\in(a,b)$.
\end{theorem}

%\textcolor{blue}{Is there any general expression to the term $\int_{x^*}^bw(x)e(x) dx$?}

\begin{proof}
The proof for the first term is exactly the same as in Theorem 10.1 of \cite{book1}. The term $\int_{x^*}^bw(x)e(x) dx$ comes from the approximation done in (\ref{fmmi2}) and it is clear that is zero, as $e(x)=0$, if $f$ is a piece-wise polynomial of degree $2l+1$.
\end{proof}
%\textcolor{blue}{For clarity, we can prove the theorem in similar lines as given in the textbook and present it in the appendix!}
%It can be seen that the number of terms of $C(x)$ in (\ref{C}), can be chosen depending on the accuracy of the quadrature rule chosen to approximate
%It is known that in Gaussian quadrature, the formulas are developed for the interval $[-1,1]$, so the theory presented in Theorem \ref{teo} can be easily adapted, just setting $[a,b]=[-1,1]$ as in Figure \ref{disc_alpha_corner_gauss}. For any other interval, the classical change of variables $t=\frac{2x-a-b}{}$

\begin{remark}
As mentioned in the previous section, one interesting case is when $w(x)=1$, as the correction can be given explicitly.\\
%\newline
%\textcolor{blue}{Can we explore this formulation for the well-known weight functions $w(x)$?\\\\
%(Look at section:other forms, at
%$https://en.wikipedia.org/wiki/Gaussian\_quadrature$}
\end{remark}

\section{Discussion about the accuracy needed when $x^*$ is approximated}\label{acc_x*}

In the case that $x^*$ is approximated, we can use Taylor expansion on (\ref{tayl_T}) to obtain some insight about which accuracy is needed in the location of the discontinuity. For simplicity, let us use here $T_{x^*}^l$ instead of $T_{x^*}^l(f^+-f^-)$. Let us consider also that we have located the position of the discontinuity but a small error of location $\beta$ has been introduced in the process,
\begin{equation}\label{erx*}
\tilde{x}^*=x^*+\beta,
\end{equation}
so that $\tilde{x}^*\in(a,b)$, being $[a,b]$ the interval of integration.
Using Taylor expansion on the jumps of the function and its derivatives in the expression that we get from (\ref{tayl_T}), and using $\tilde{x}^*$ instead of $x^*$, 

\begin{equation}\label{tayl_T_tilde}
T_{\tilde{x}^*}^l(f^+-f^-)(x)=\sum_{k=0}^{l}\frac{\left[f^{(k)}(\tilde{x}^*)\right]}{k!}(x-\tilde{x}^*)^k,
\end{equation}
we can state the following Theorem:

\begin{theorem}\label{teop}
    The precision needed in the location of the discontinuity must be $O(b-x^*)$.
\end{theorem}
\begin{proof}
    If we assume that there is an error in the location of $x^*$, as expressed in (\ref{erx*}), we can use Taylor expansions in the jumps of the function in (\ref{tayl_T}),
    $$[f^{(l)}(\tilde{x}^*)]=\sum_{k=0}^{m}\frac{\left[f^{(l+k)}(x^*)\right]}{k!}\beta^k.$$
    Replacing in (\ref{tayl_T_tilde}) and simplifying, we obtain that,
\begin{equation}\label{res_teo}
    T_{\tilde{x}^*}^l(f^+-f^-)(x)=T_{x^*}^l(f^+-f^-)(x)+O(\beta(x-x^*)^l).
\end{equation}
As the correction is applied in the interval $[x^*,b]$, then it should be clear that the error in the location of $x^*$ implies an error for $T_{x^*}^l(f^+-f^-)(x)$ in (\ref{res_teo}) that is $O(\beta (b-x^*)^l)$. From this last expression, if $\beta=O(b-x^*)$ we obtain the accuracy needed for $T_{\tilde{x}^*}^l(f^+-f^-)(x)$.
    
\end{proof}

\begin{remark}
If the jumps in the function and its derivatives at $x^*$ or $\tilde{x}^*$ (if Theorem \ref{teop} is satisfied) are approximated, then it is straightforward to see from the expressions in (\ref{res_teo}) and (\ref{tayl_T}) that the accuracy needed if the jumps in the function and the derivatives  $[f^{(k)}(x^*)], k=0,\cdots, l$ are approximated, is $O((b-x^*)^{l+1-k})$ for $k=0,\cdots, l$.
\end{remark}

\section{Numerical experiments}\label{num_exp}
In this section we will present some experiments to check the application of the  correction terms introduced in the previous sections to some particular integration rules.

%\textcolor{blue}{We have to try some numerical results for the proposed theory!}
\subsection{Experiments for Newton-Cotes formulas}
In the case of Newton-Cotes integration rules, the data is typically provided at discrete points, often distributed over a uniform grid. To express the data on the $+$ side in terms of the data on the $-$ side (or vice versa), it becomes necessary to use formula (\ref{fmmi1}) prior to computing the numerical integration. When we use the formula in (\ref{fmmi1}), the resulting function exhibits smoothness up to the degree $l$ specified in the formula. This ensures that all the derivatives up to the $l^{th}$ term will be continuous, providing a smooth representation of the function. Following this step, we apply the classical numerical integration rule and we add the correction term. This process is crucial as it enables us to accurately approximate the integral, even when the data is only available at isolated points. Let us check now that this approach enhances the precision of the Newton-Cotes integration rules. In what follows, we use a grid refinement analysis to check the numerical accuracy.

In the grid refinement experiments, we compute the error $E_i$ as the absolute difference between the exact integral and its approximation, which is derived using composite quadrature rules. The accuracy order is determined by the formula,
\begin{equation}\label{orden}
O_i=\frac{\ln(E_i/ E_{i+1})}{\ln(n_i/n_{i+1})},
\end{equation}
where $E_i$ represents the error for a grid with $n_i$ points and $E_{i+1}$ denotes the error for a grid with $n_{i+1}$ points. %In these experiments, for the trapezoid rule and the Simpson's $\frac{3}{8}$ rule, $n_{i+1}$ is twice $n_i$, while for the Simpson's $\frac{1}{3}$ rule, $n_{i+1}$ equals $2n_i+1$.
\begin{figure}
\begin{center}
\includegraphics[height=4.5cm]{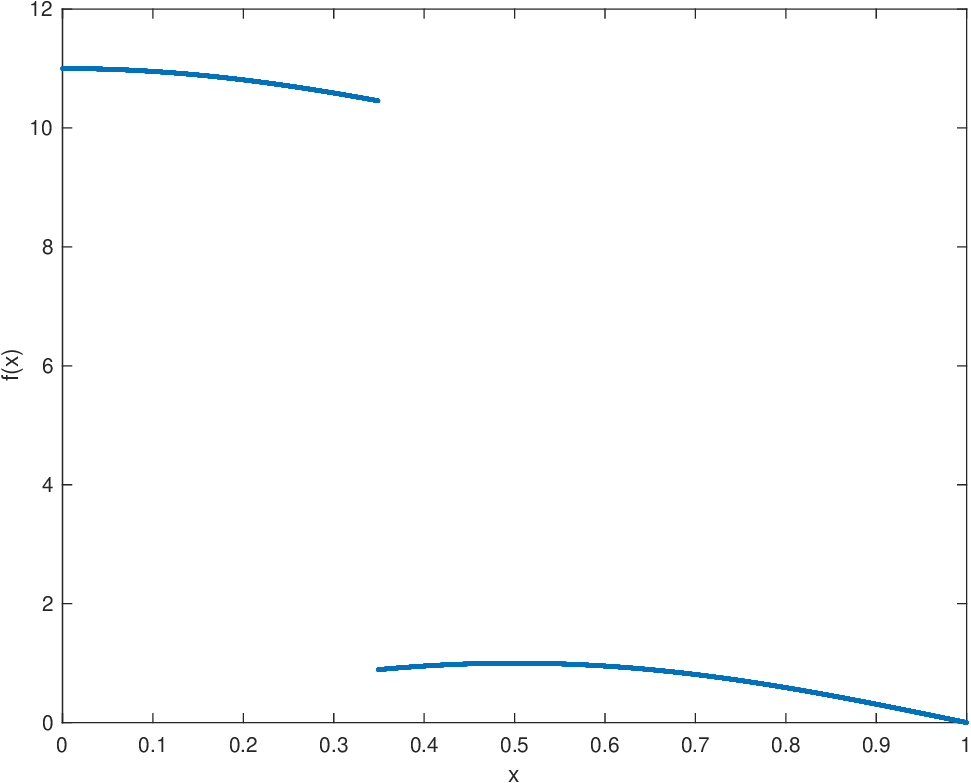}
\includegraphics[height=4.5cm]{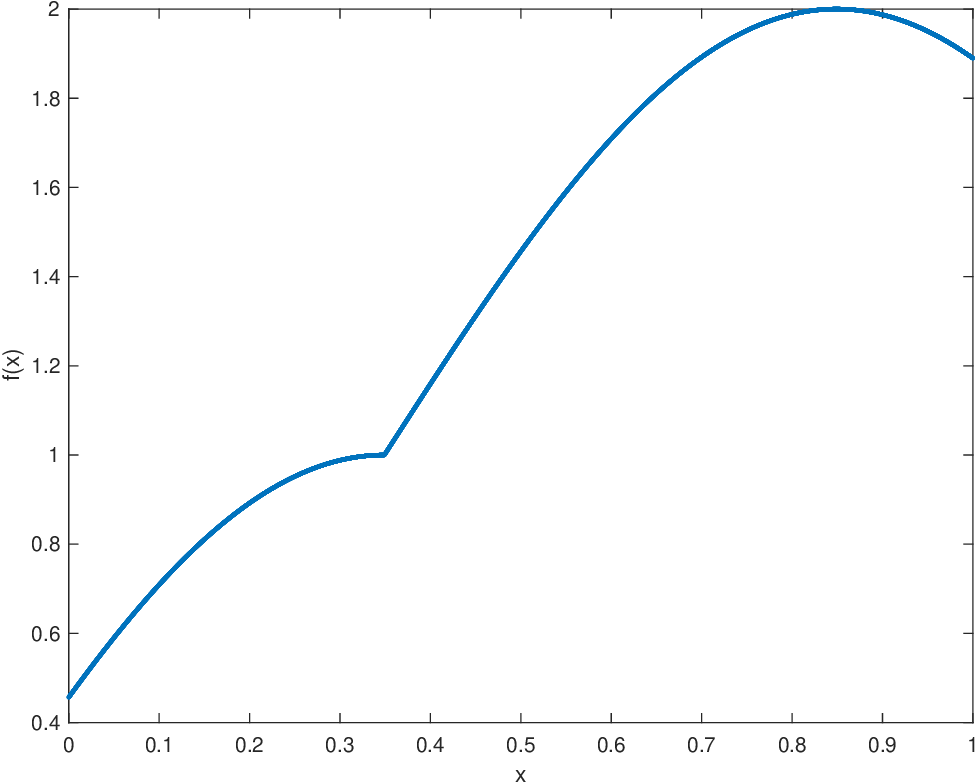}
\includegraphics[height=4.5cm]{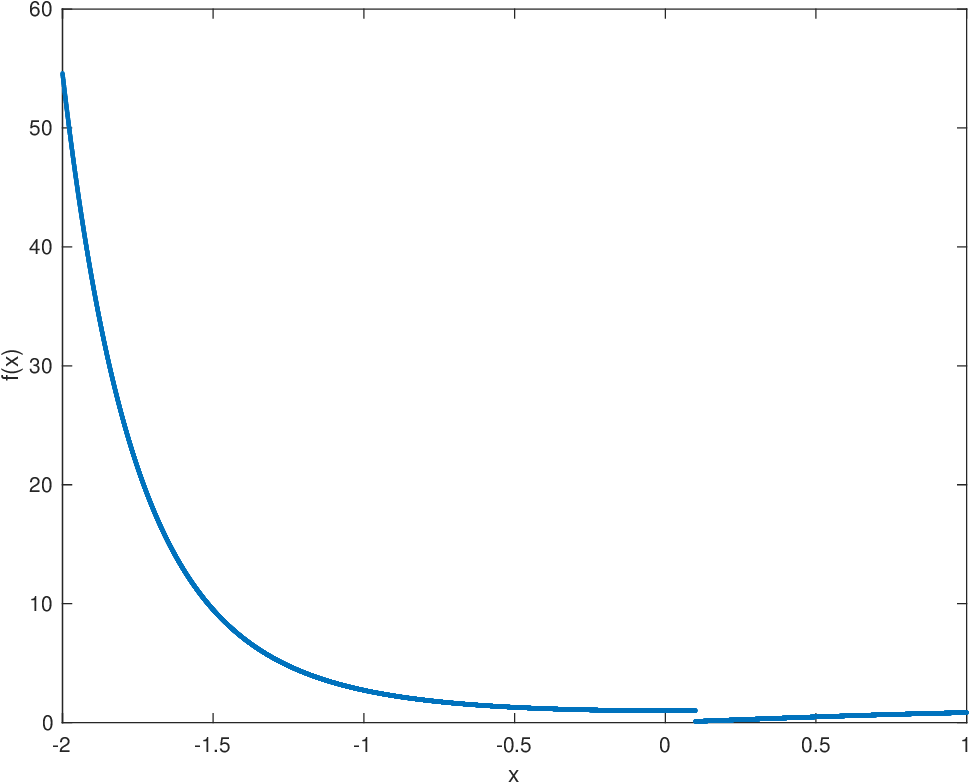}
\caption{Functions presented in  (\ref{function}), (\ref{function_corner}), and (\ref{func_comp})used in Experiments 1, 2, and 4.}\label{funciones}
\end{center}
\end{figure}

\subsubsection{Experiment 1}
Let us consider the following piece-wise continuous function
\begin{equation}\label{function}
%\scalemath{0.9}{
f(x)=\left\{\begin{array}{ll}
\cos\left(\pi x\right) + 10, & \textrm{ if } 0\le x<\frac{\pi}{9},\\
\sin\left(\pi x\right), & \textrm{ if } \frac{\pi}{9}\le x\le 1.
\end{array}
\right.
%}
\end{equation}
This function has been plotted to the left in Figure \ref{funciones}.

\begin{figure}
\begin{center}
%quadrature_IIM/programas/refinamientos_Newton_Cotes_nuevo_articulo.m
%quadrature_IIM/programas/simpson_articulo_R1.m

\includegraphics[height=4.5cm]{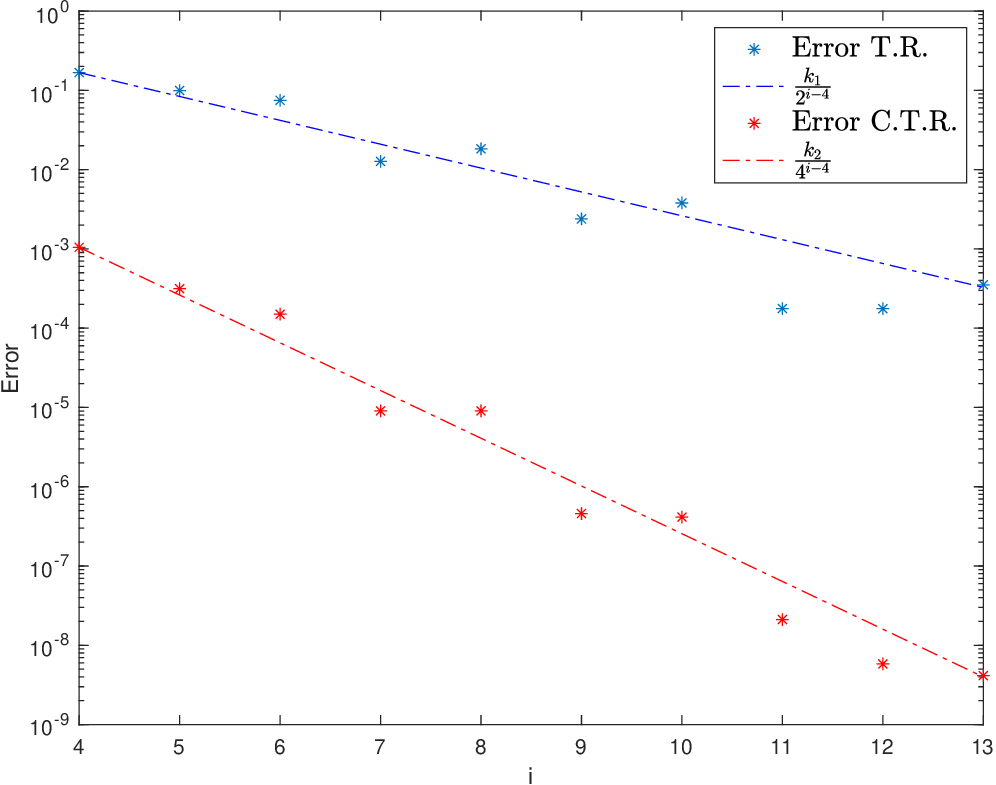}
\includegraphics[height=4.5cm]{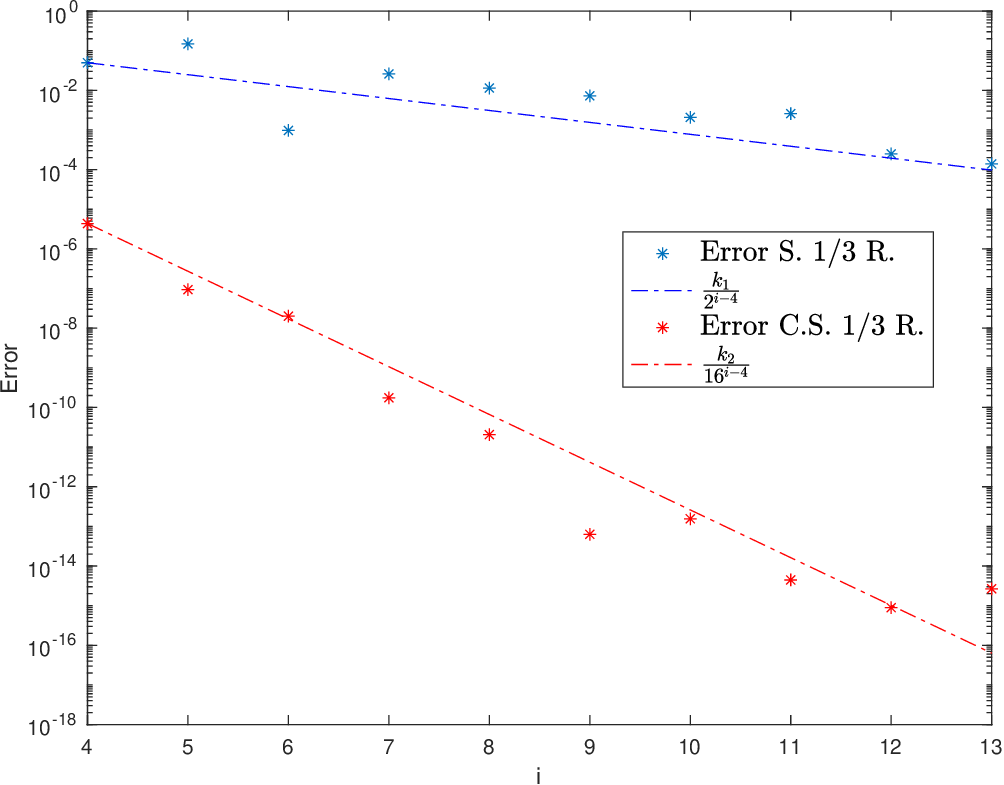}
\includegraphics[height=4.5cm]{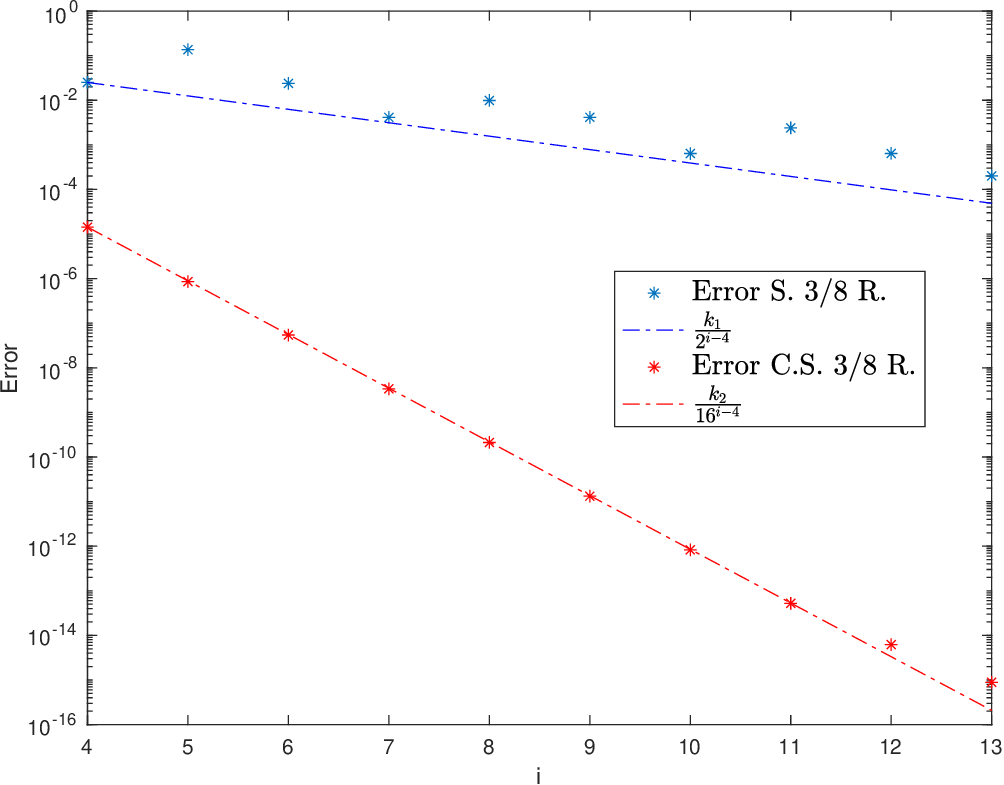}
\caption{Grid refinement analysis for the numerical integration of the function in (\ref{function}). To the left, using the composite trapezoid rule and the corrected composite trapezoid rule. At the center, using the composite Simpson's $1/3$ rule and the corrected rule. To the right, using the composite Simpson's $3/8$ rule and the corrected rule. We can see that the theoretical error decreases with $i$, being $n=2^i$ the number of nodes used. The error of the corrected formulas decreases following the theoretical rate. Tables \ref{conv}, \ref{s_1_3} and \ref{s3_8} contain the data represented in these graphs.}\label{plot_err}
\end{center}
\end{figure}

In Figure \ref{plot_err} we show a grid refinement analysis for the integration of the function in (\ref{function}) using the classical trapezoid rule, the Simpson's 1/3 and the Simpson's 3/8 rules, jointly with their corrected versions proposed in Section \ref{corr_terms_general}. Tables \ref{conv}, \ref{s_1_3} and \ref{s3_8} contain the data represented in these graphs. We can see how the error of the corrected trapezoid rule decreases with $O(h^2)$ accuracy, and the error for Simpson's 1/3 and 3/8 rules decrease with $O(h^4)$ accuracy. The classical formulas attain $O(h)$ accuracy due to the presence of the jump discontinuity of $f(x)$ at $x=\frac{\pi}{9}$. The dashed lines in the graphs represent the expected theoretical slope for the decreasing of the errors attained by each rule.

%\begin{landscape}
\begin{table}[!ht]
\begin{center}
\resizebox{18cm}{!} {
\begin{tabular}{|c|c|c|c|c|c|c|c|c|c|c|c|c|c|c|c|c|c|c|c|c|c|c|c|c|c|c|}
%\hline $x=i-4$ & 0  & 1      & 2     & 3    & 4 & 5  & 6 &7 &8&9%&$2^{14}$&$2^{15}$
%\\
\hline $n=2^i$ & $2^4$  & $2^5$      & $2^6$     & $2^7$    & $2^8$ & $2^9$  & $2^{10}$ &$2^{11}$ &$2^{12}$&$2^{13}$%&$2^{14}$&$2^{15}$
      \\       
\hline
Error T.R. ($E_i^{\infty}$)& 1.67339e-01 & 9.93916e-02 & 7.47340e-02 & 1.26944e-02 & 1.83072e-02 & 2.38379e-03 & 3.79302e-03 & 1.76665e-04 & 1.76013e-04 & 3.52282e-04 
      \\       
\hline 
$O_i$&- & 0.75158 & 0.41136 & 2.5576 & -0.52822 & 2.9411 & -0.67009 & 4.4243 & 0.0053312 & -1.001 
    \\    
\hline     
Error C.T.R. ($E_i^{\infty}$)& 1.04638e-03 & 3.14790e-04 & 1.50082e-04 & 9.04721e-06 & 9.06771e-06 & 4.59216e-07 & 4.14277e-07 & 2.11422e-08 & 5.81906e-09 & 4.14297e-09
      \\       
\hline     
$O_i$& -& 1.7329 & 1.0686 & 4.0521 & -0.0032652 & 4.3035 & 0.14858 & 4.2924 & 1.8613 & 0.49012
\\
\hline    
\end{tabular}
}
\caption{Grid refinement analysis in the infinity norm for the integral of the function shown in (\ref{function}) using the composite trapezoid rule. The table shows the results for the trapezoidal rule (T.R.) and the corrected trapezoidal rule (C.T.R.). %The bottom part shows the Simpson's 1/3 Rule (S. 1/3 R.) and the corrected Simpson's rule (C.S. 1/3 R.).
}\label{conv}
\end{center}
\end{table}

%\begin{landscape}
\begin{table}[!ht]
\begin{center}
\resizebox{18cm}{!} {
\begin{tabular}{|c|c|c|c|c|c|c|c|c|c|c|c|c|c|c|c|c|c|c|c|c|c|c|c|c|c|c|}
\hline $n=2^i$ & $2^4$  & $2^5$      & $2^6$     & $2^7$    & $2^8$ & $2^9$  & $2^{10}$ &$2^{11}$ &$2^{12}$&$2^{13}$%&$2^{14}$&$2^{15}$
      \\       
\hline
Error S. 1/3 R. ($E_i^{\infty}$)& 4.96488e-02 & 1.48172e-01 & 9.71040e-04 & 2.59437e-02 & 1.14289e-02 & 7.25740e-03 & 2.08575e-03 & 2.58583e-03 & 2.50037e-04 & 1.39233e-04
      \\       
\hline 
$O_i$&- & -1.5774 & 7.2535 & -4.7397 & 1.1827 & 0.65516 & 1.7989 & -0.31006 & 3.3704 & 0.84464
    \\
\hline     
Error C.S. 1/3 R. ($E_i^{\infty}$)& 4.36231e-06 & 9.35406e-08 & 2.01092e-08 & 1.75421e-10 & 2.06226e-11 & 6.30607e-14 & 1.55431e-13 & 4.44089e-15 & 8.88178e-16 & 2.66454e-15
      \\       
\hline     
$O_i$& -& 5.5434 & 2.2177 & 6.8409 & 3.0885 & 8.3533 & -1.3015 & 5.1293 & 2.3219 & -1.585
\\
\hline    
\end{tabular}
}
\caption{Grid refinement analysis in the infinity norm for the integral of the function shown in (\ref{function}) using the composite Simpson 1/3 rule. The table shows the results for the Simpson 1/3 rule (S. 1/3 R.) and the corrected Simpson 1/3 rule (C.S. 1/3 R.).
}\label{s_1_3}
\end{center}
\end{table}

%\begin{landscape}
\begin{table}[!ht]
\begin{center}
\resizebox{18cm}{!} {
\begin{tabular}{|c|c|c|c|c|c|c|c|c|c|c|c|c|c|c|c|c|c|c|c|c|c|c|c|c|c|c|}
\hline $n=2^i$ & $2^4$  & $2^5$      & $2^6$     & $2^7$    & $2^8$ & $2^9$  & $2^{10}$ &$2^{11}$ &$2^{12}$&$2^{13}$%&$2^{14}$&$2^{15}$
      \\       
\hline
Error S. 3/8 R. ($E_i^{\infty}$)& 2.50247e-02 & 1.35857e-01 & 2.38721e-02 & 4.12512e-03 & 9.87350e-03 & 4.14324e-03 & 6.39055e-04 & 2.39115e-03 & 6.39339e-04 & 2.01386e-04
      \\       
\hline 
$O_i$&- & -2.4407 & 2.5087 & 2.5328 & -1.2591 & 1.2528 & 2.6967 & -1.9037 & 1.9031 & 1.6666
    \\
\hline     
Error C.S. 3/8 R. ($E_i^{\infty}$)& 6.10369e-06 & 3.83753e-07 & 2.41635e-08 & 1.50811e-09 & 9.40217e-11 & 5.88241e-12 & 3.66818e-13 & 2.39808e-14 & 2.66454e-15 & 8.88178e-16
      \\       
\hline     
$O_i$& -& 3.9914 & 3.9893 & 4.002 & 4.0036 & 3.9985 & 4.0033 & 3.9351 & 3.1699 & 1.585
\\
\hline    
\end{tabular}
}
\caption{Grid refinement analysis in the infinity norm for the integral of the function shown in (\ref{function}) using the composite Simpson 3/8 rule. The table shows the results for the Simpson 3/8 rule (S. 3/8 R.) and the corrected Simpson 3/8 rule (C.S. 3/8 R.).
}\label{s3_8}
\end{center}
\end{table}

\subsubsection{Experiment 2}

Let us now consider an example of a continuous function with jumps in the derivatives:

%y1=cos(pi*(x1-pos_disc))+a;
%y2=1+sin(pi*(x2-pos_disc));
\begin{equation}\label{function_corner}
%\scalemath{0.9}{
f(x)=\left\{\begin{array}{ll}
\cos\left(\pi (x-\frac{\pi}{9})\right), & \textrm{ if } 0\le x<\frac{\pi}{9},\\
1+\sin\left(\pi (x-\frac{\pi}{9})\right), & \textrm{ if } \frac{\pi}{9}\le x\le 1.
\end{array}
\right.
%}
\end{equation}
This function has been plotted at the center of Figure \ref{funciones}. Our aim is to show that, in this case, we can locate the position of the discontinuity and approximate the jump in the function and the derivatives assuring that the accuracy of the correction is kept.
In this case, we use the algorithm in \cite{Cohen-arandiga-donat-dyn} to locate the discontinuity, and we have approximated the jumps in the function and its derivatives at the approximated location of the discontinuity using one-sided interpolation, (see for example \cite{IIM_newton}). For this continuous function, the trapezoid rule already attains a good level of accuracy, as a jump in the first derivative of the function turns into an error of $O(h^2)$ for its integral, that is the accuracy for the composite trapezoidal rule. Thus, we present the results for the Simpson's $1/3$ rule.

\begin{figure}
\begin{center}
\includegraphics[height=4.5cm]{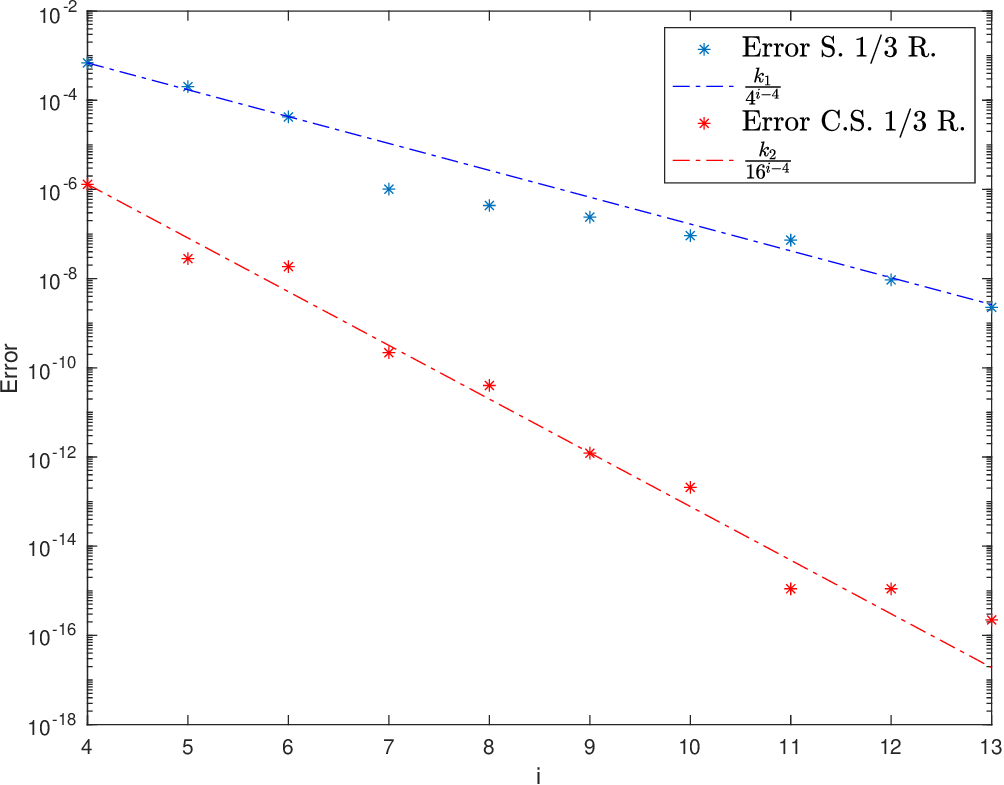}
\caption{Grid refinement analysis for the numerical integration of the function in (\ref{function_corner}). In this example we use the composite Simpson's $1/3$ rule and the corrected rule. We can see that the error of the corrected formulas decrease following the theoretical rate.}\label{plot_err_corner}
\end{center}
\end{figure}

In Figure \ref{plot_err_corner} we present a grid refinement analysis for the integration of the function in (\ref{function_corner}) using the classical Simpson's $1/3$ rule and the general correction proposed in Section \ref{corr_terms_general}. Table \ref{s_1_3_corner} presents the results obtained. Again, we can see that the results show a decreasing of the error with slopes that are similar to the expected theoretical ones, and that are represented through  dashed lines in Figure \ref{plot_err_corner}. We can see that the classical Simpson's $1/3$ rule decreases with $O(h^2)$ accuracy in this case, while the corrected rule decreases with $O(h^4)$.

%\begin{landscape}
\begin{table}[!ht]
%quadrature_IIM/programas/simpson_articulo_R1_corner.m
\begin{center}
\resizebox{18cm}{!} {
\begin{tabular}{|c|c|c|c|c|c|c|c|c|c|c|c|c|c|c|c|c|c|c|c|c|c|c|c|c|c|c|}
\hline $n=2^i$ & $2^4$  & $2^5$      & $2^6$     & $2^7$    & $2^8$ & $2^9$  & $2^{10}$ &$2^{11}$ &$2^{12}$&$2^{13}$%&$2^{14}$&$2^{15}$
      \\       
\hline
Error S. 1/3 R. ($E_i^{\infty}$)& 6.86434e-04 & 2.01877e-04 & 4.21341e-05 & 1.01332e-06 & 4.36146e-07 & 2.39839e-07 & 9.20051e-08 & 7.31517e-08 & 9.33086e-09 & 2.26839e-09
      \\       
\hline 
$O_i$&- &1.7656 & 2.2604 & 5.3778 & 1.2162 & 0.86275 & 1.3823 & 0.33082 & 2.9708 & 2.0403
    \\
\hline     
Error C.S. 1/3 R. ($E_i^{\infty}$)&1.30320e-06 & 2.80621e-08 & 1.85818e-08 & 2.18269e-10 & 4.01432e-11 & 1.21636e-12 & 2.08500e-13 & 1.11022e-15 & 1.11022e-15 & 2.22045e-16
      \\       
\hline     
$O_i$& -& 5.5373 & 0.59473 & 6.4116 & 2.4429 & 5.0445 & 2.5445 & 7.5531 & 0 & 2.3219
\\
\hline    
\end{tabular}
}
\caption{Grid refinement analysis in the infinity norm for the integral of the function shown in (\ref{function_corner}) using the composite Simpson 1/3 rule. The table shows the results for the Simpson 1/3 rule (S. 1/3 R.) and the corrected Simpson 1/3 rule (C.S. 1/3 R.).
}\label{s_1_3_corner}
\end{center}
\end{table}

\subsection{Experiments for Gauss-Legendre quadrature formulas}

In this section we check if the new technique allows for the recovering of the reproduction of polynomials for Gauss-Legendre quadrature formulas when the data is piece-wise smooth. 

\subsubsection{Experiment 3}
In this third experiment, we check the Gauss-Legendre integration rules over the interval \([-1, 1]\), given by the expression
\[
\int_{-1}^1 f(x) \, dx \approx \sum_{i=1}^n w_i f(x_i),
\]
where \(x_i\) are the roots of the \(n\)-th Legendre polynomial \(P_n(x)\), and the weights \(w_i\) are calculated as:
\[
w_i = \frac{2}{(1 - x_i^2) [P_n'(x_i)]^2}.
\]
We know that this formula ensures that the approximation is exact for polynomials of degree up to \(2n-1\). We consider that the integrand is a piece-wise polynomial of the degree that the formula reproduces. Our aim is to check if, introducing a correction term, we can recover the accuracy of the original formula. Table \ref{tabla_GL} presents the classical Gauss-Legendre points $x_i$ and weights $w_i$ that we use in the experiment.

\begin{table}[!ht]
\begin{center}
\resizebox{8cm}{!} {
\begin{tabular}{|c|c|c|}
\hline
Number of points, \(n\) & Points, \(x_i\) & Weights, \(w_i\) \\
\hline
%1 & 0 & 2 \\
%\hline
2 & \(\pm \frac{1}{\sqrt{3}}\) & 1 \\
\hline
3 & 0 & \(\frac{8}{9}\) \\
  & \(\pm \sqrt{\frac{3}{5}}\) & \(\frac{5}{9}\) \\
\hline
4 & \(\pm \sqrt{\frac{3}{7} - \frac{2}{7} \sqrt{\frac{6}{5}}}\) & \(\frac{18 + \sqrt{30}}{36}\) \\
  & \(\pm \sqrt{\frac{3}{7} + \frac{2}{7} \sqrt{\frac{6}{5}}}\) & \(\frac{18 - \sqrt{30}}{36}\) \\
\hline
5 & 0 & \(\frac{128}{225}\) \\
  & \(\pm \frac{1}{3} \sqrt{5 - 2 \sqrt{\frac{10}{7}}}\) & \(\frac{322 + 13\sqrt{70}}{900}\) \\
  & \(\pm \frac{1}{3} \sqrt{5 + 2 \sqrt{\frac{10}{7}}}\) & \(\frac{322 - 13\sqrt{70}}{900}\) \\
\hline
\end{tabular}
}
\caption{ Classical Gauss-Legendre points $x_i$ and weights $w_i$.}\label{tabla_GL}
\end{center}
\end{table}

\begin{figure}
\begin{center}
\includegraphics[height=5cm]{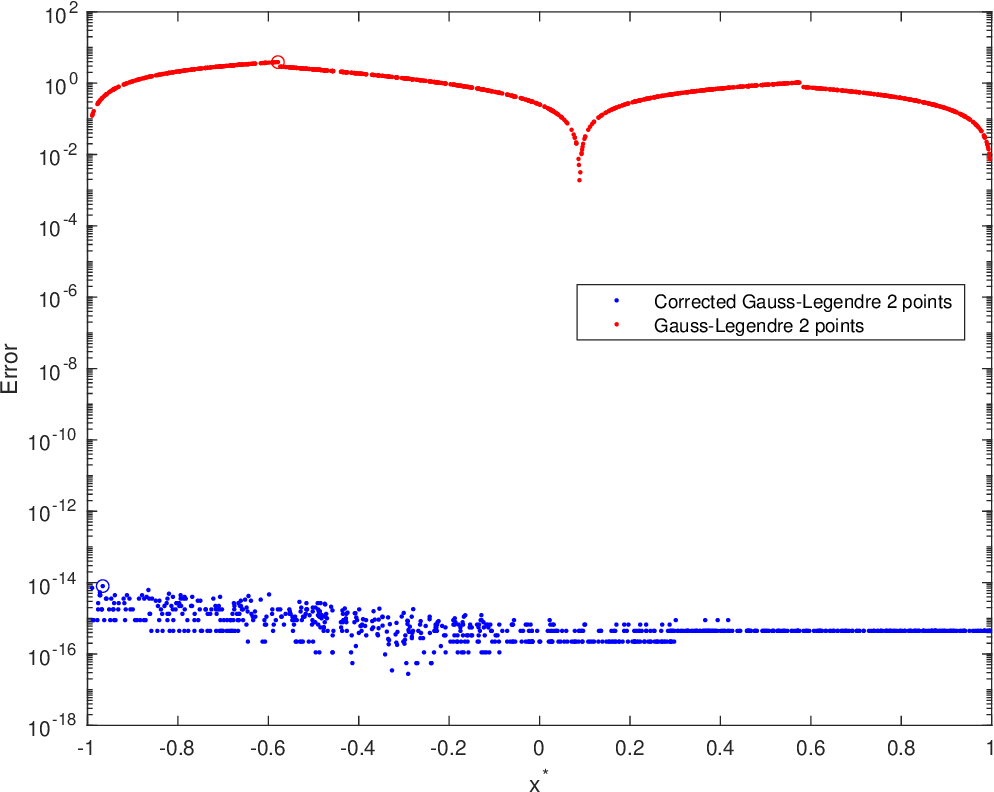}
\includegraphics[height=5cm]{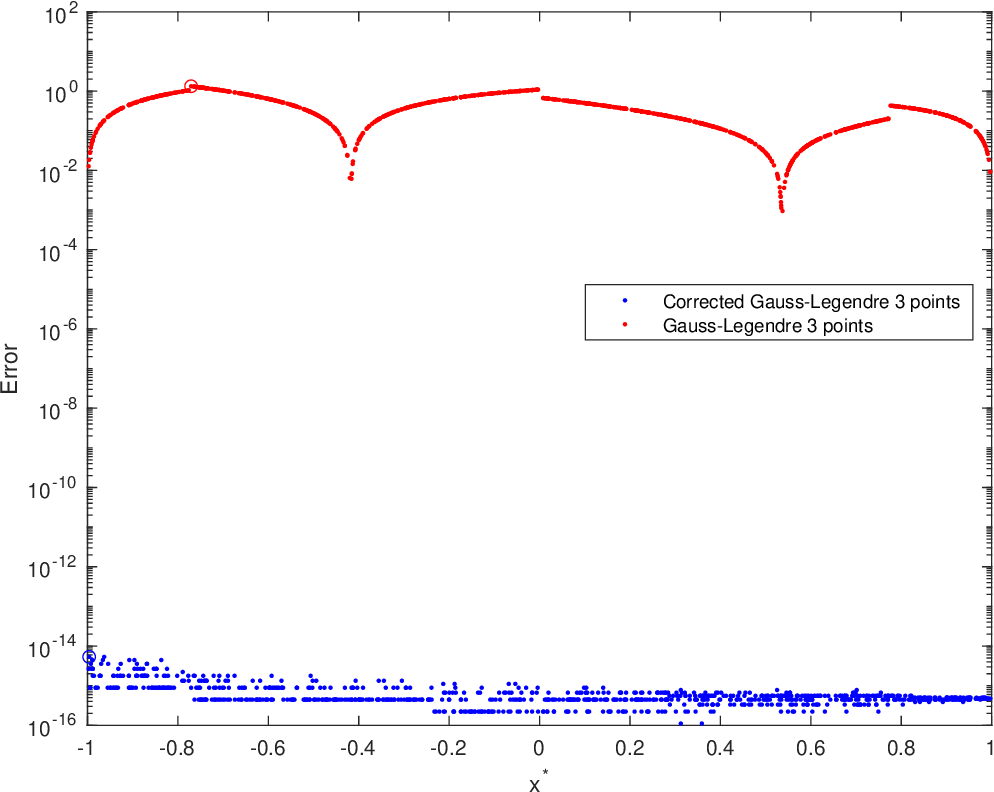}\\
\includegraphics[height=5cm]{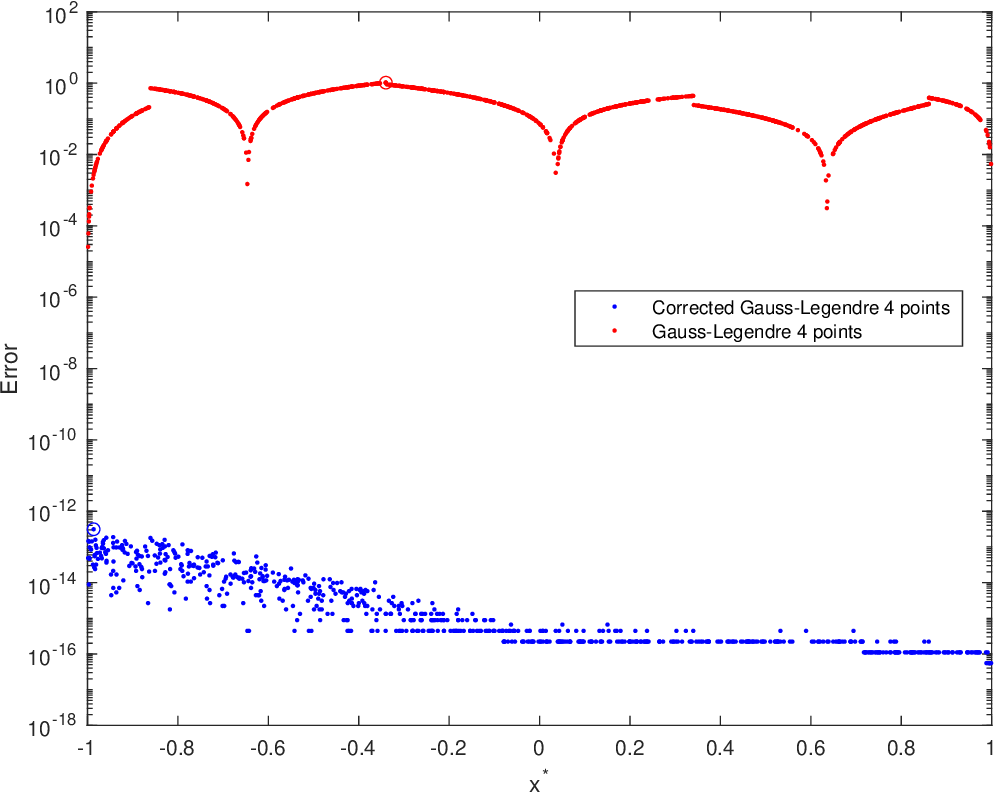}
\includegraphics[height=5cm]{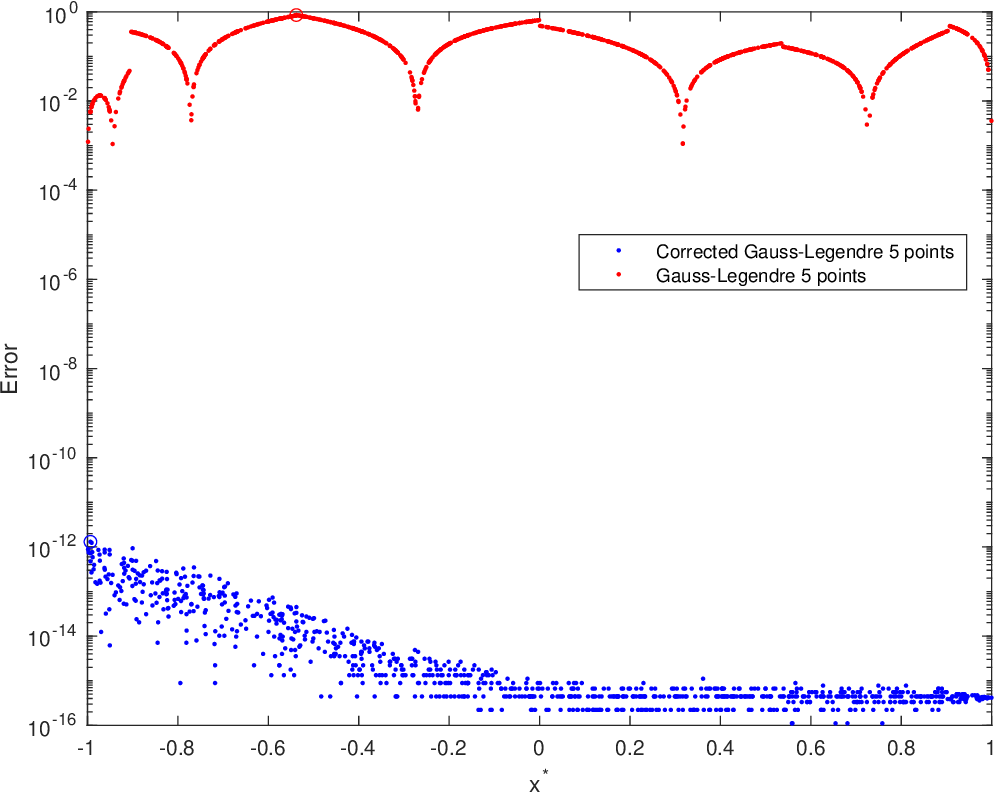}
\caption{Representation of the error obtained when choosing randomly 1000 samples for $x^*\in[-1,1]$ and integrating the piece-wise polynomial functions in Table (\ref{polinomios}) using the Gaussian quadrature formulas in Table (\ref{tabla_GL}), with and without correction terms. In the $x$ axis we have represented the position of the discontinuity $x^*$ and in the $y$ axis, the error obtained with the corresponding Gauss-Legendre quadrature formula. With a circle we have marked the maximum error obtained for the 1000 experiments performed for the corrected and classical formula.}\label{err_GL}
\end{center}
\end{figure}

%with
%\begin{equation}\label{G1}
%A_0=1, \quad A_1=1, \quad x_1=-\frac{\sqrt{3}}{3},\quad x_2=\frac{\sqrt{3}}%{3},
%\end{equation}
%Let us consider the formulas in Table \ref{tabla_GL}. 
In what follows, we check that, including the correction term explained in Section \ref{corr_terms_general}, we obtain an exact result up to the machine precision for the piece-wise polynomials shown in Table \ref{polinomios},
with $x^*$ being a random number in the interval of integration $[-1,1]$. The experiment has been set as we describe in what follows: We perform 1000 experiments for each integration rule shown in Table \ref{tabla_GL}, where a random position of the discontinuity $x^*$ is chosen for each experiment. The error for each different $x^*$, that has been represented in the $x$ axis,  is plotted in Figure \ref{err_GL} for each quadrature rule. Then we choose the maxima of the error for the 1000 experiments for each rule, that have been represented in Figure \ref{err_GL} with a blue circle for the corrected formula and with a red circle for the classical ones, and we present these errors in Table \ref{tab_exp1}. We can see that the maxima of the error in all the cases is close to the machine precision for the corrected formulas.

\begin{table}[ht]
\centering
\begin{tabular}{|c| >{\centering\arraybackslash}m{12cm}|}
\hline
Number of points, \( n \) & Piecewise polynomial \( f(x) \) \\
\hline
2 & 
\(
f(x)=\left\{\begin{array}{ll}
 x^3 + 2x^2 - 3x + 1, & x < x^* \\
 2x^3 - 2x^2 + x - 2, & x \geq x^*
\end{array}\right.
\) \\
\hline
3 & 
\(
f(x)=\left\{\begin{array}{ll}
 x^5 - 3x^4 + x^3 - x^2 + x + 1, & x < x^* \\
 2x^5 - x^4 + 2x^3 - x^2 - 2x + 3, & x \geq x^*
\end{array}\right.
\) \\
\hline
4 & 
\(
f(x)=\left\{\begin{array}{ll}
 -x^7 + x^6 + x^5 - 3x^4 + x^3 - x^2 + x + 1, & x < x^* \\
 2x^7 - x^6 + 2x^5 - x^4 + 2x^3 - x^2 - 2x + 3, & x \geq x^*
\end{array}\right.
\) \\
\hline
5 & 
\(
f(x)=\left\{\begin{array}{ll}
 x^9-2x^8-x^7 + x^6 + x^5 - 3x^4 + x^3 - x^2 + x + 1, & x < x^* \\
 3x^9-x^8+2x^7 - x^6 + 2x^5 - x^4 + 2x^3 - x^2 - 2x + 3, & x \geq x^*
\end{array}\right.
\) \\
\hline
\end{tabular}\caption{Piece-wise polynomials defined for the different number of points $n$ considered.}\label{polinomios}
\end{table}

%\begin{equation}\label{exp1}
%f(x)=\left\{\begin{array}{cc}
% x^3+2x^2-3x+1,    & x<x^*, \\
% 2x^3-2x^2+x-2,    & x\geq x^*,
%\end{array}\right.
%\end{equation}

\begin{table}[!ht]
\begin{center}
%\resizebox{18cm}{!} {
\begin{tabular}{|c|c|c|}
\hline $n$ & Error Classical & Error Corrected  \\
\hline
     2 & 3.8906 & 7.9936e-15\\
\hline
     3 &1.3243 & 5.3291e-15\\
\hline
     4 &1.0313 & 3.1353e-13\\
\hline
     5 &0.84541 & 1.3056e-12\\  
\hline
\end{tabular}
%}
\caption{Maximum of the error obtained when choosing randomly 1000 samples for $x^*\in[-1,1]$ and integrating the piece-wise polynomial functions in Table (\ref{polinomios}) using the Gaussian quadrature formulas in Table (\ref{tabla_GL}), with and without correction terms.
}\label{tab_exp1}
\end{center}
\end{table}

\subsubsection{Experiment 4}

In this subsection we want to test the performance of composed Gauss-Legendre quadrature rules. The idea is to obtain an approximation of an integral over an arbitrary interval using a uniform grid and a quadrature rule for each subinterval of the grid. This can be done easily performing a change of variables for each subinterval. Let us suppose that we want to integrate the function $f(y)$ in the interval $[y_0, y_0+h]$, being $h$ the grid-spacing of a uniform grid. Then, if we use the change of variables $x=\frac{y-(y_0+h/2)}{h/2}$, we obtain
\begin{equation}\label{exp4_1}
    \int_{y_0}^{y_0+h}f(y)dy=\frac{h}{2}\int_{-1}^{1}f\left(x\frac{h}{2}+y_0+\frac{h}{2}\right)dy,
\end{equation}
that using a Gauss-Legendre quadrature rule, transforms into,
\begin{equation}\label{exp4_2}
    \int_{y_0}^{y_0+h}f(y)dy\approx\sum_{k=0}^l\frac{h}{2}W_k f(x_k\frac{h}{2}+y_0+\frac{h}{2}).
\end{equation}
Then we can use this quadrature formula on each sub-interval of the uniform grid.
The idea is to use the correction proposed in (\ref{int1}) in the sub-intervals of the composed quadrature rule affected by a discontinuity. 

\begin{figure}
\begin{center}
\includegraphics[height=5cm]{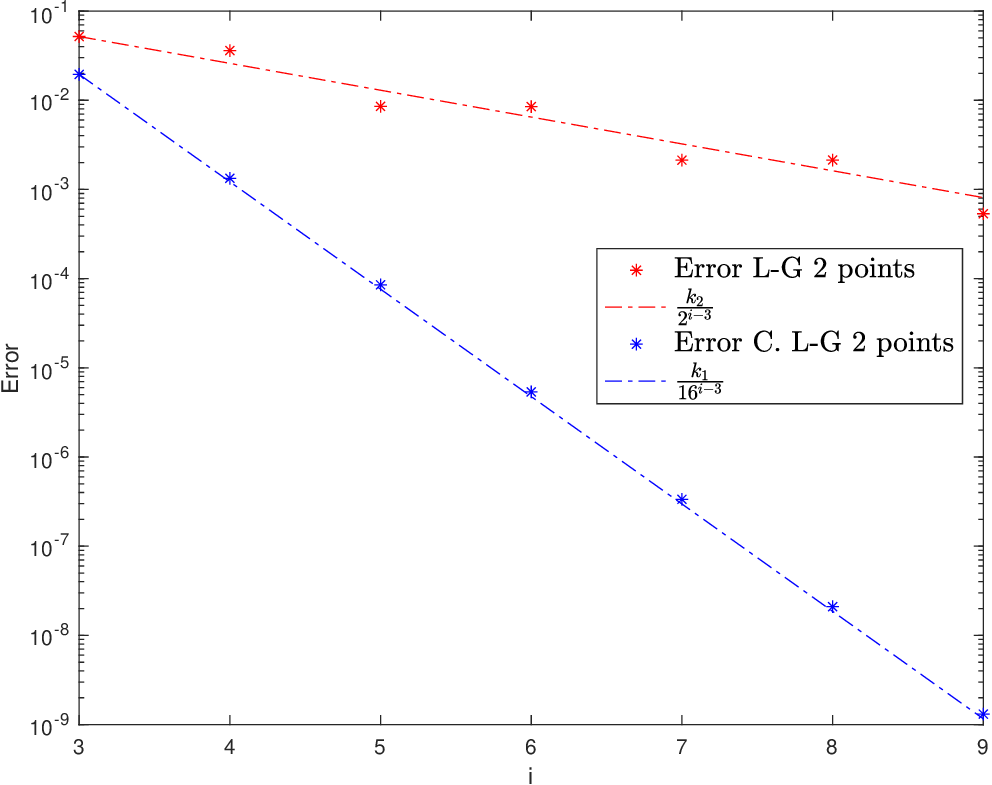}
\includegraphics[height=5cm]{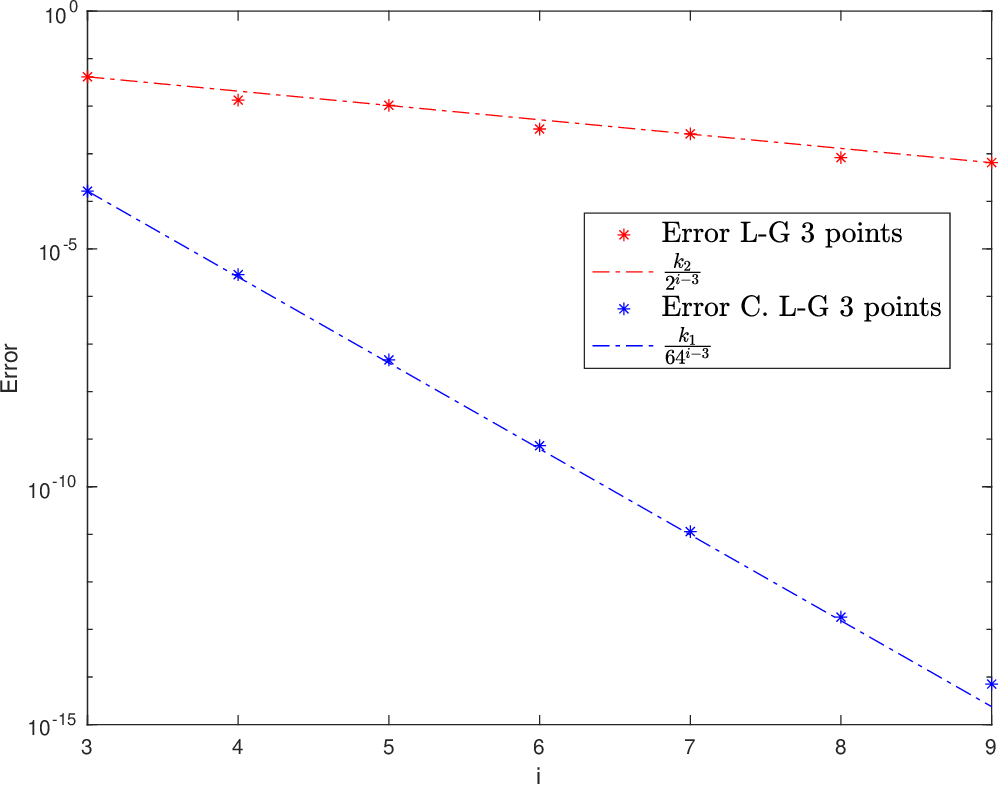}\\
\includegraphics[height=5cm]{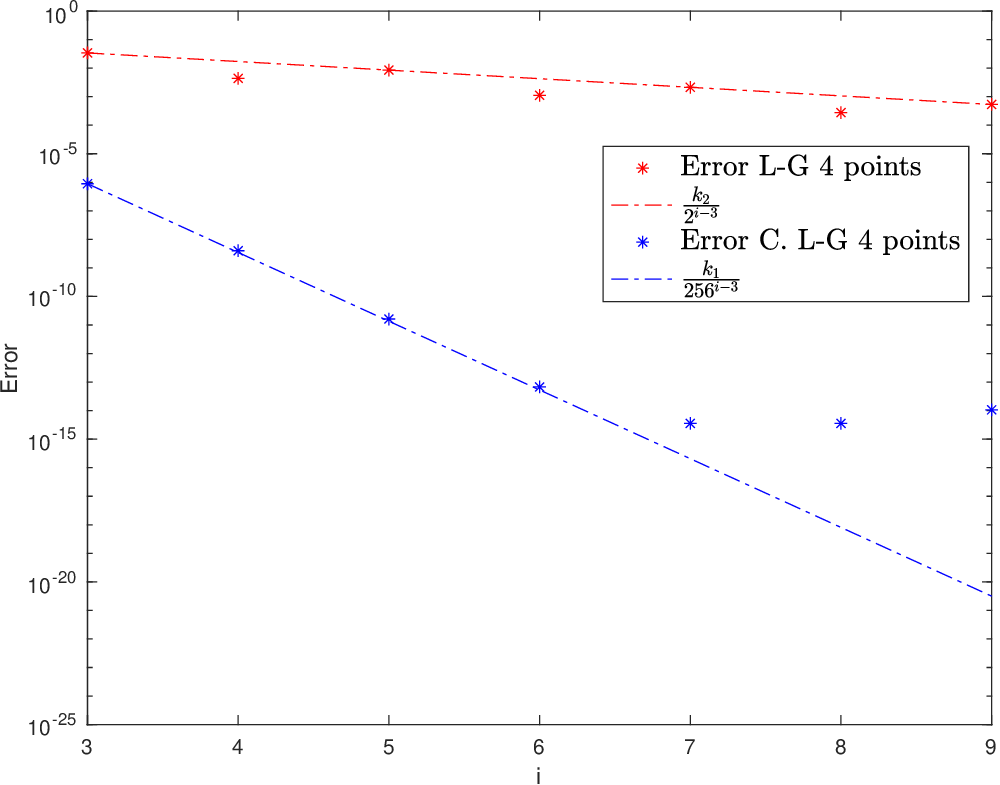}
\includegraphics[height=5cm]{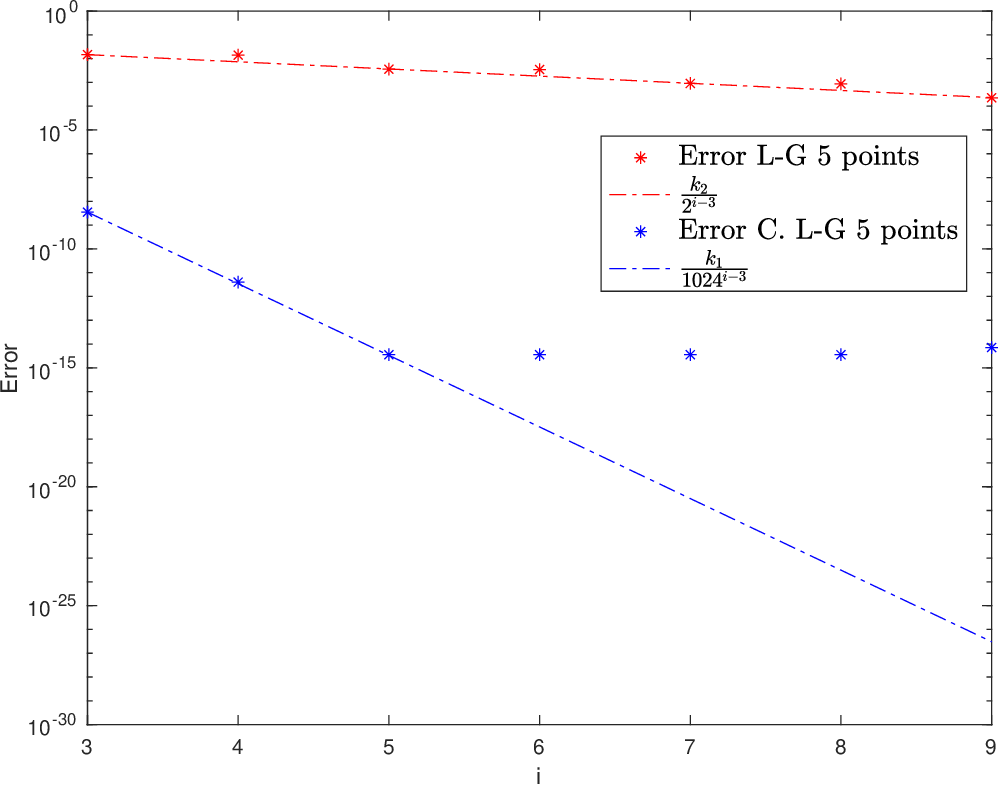}
\caption{Representation of the grid refinement analysis performed for composed Gauss-Legendre quadrature formulas for the function in (\ref{func_comp}).}\label{err_GL_composed}
\end{center}
\end{figure}

In this case we use the following piece-wise defined function,
\begin{equation}\label{func_comp}
    f(x)=\left\{\begin{array}{ll}
 e^{x^2}, & -2\le x \leq x^*, \\
 \sin(x), & x^*< x\le 1,
\end{array}\right.
\end{equation}
with $x^*=0.1$, that has been represented in Figure \ref{funciones} to the right.
We approximate the following integral
$$\int_{-2}^1 f(x)dx,$$
using Gauss-Legendre quadrature rules of different order. We perform a grid refinement analysis considering $n+1$ points in the interval $[-2,1]$, with an initial number of nodes $n=8$, and doubling $n$ in several consecutive steps. Then we obtain the absolute error of the numerical integration versus the exact integration in the interval $[-2,1]$, and obtain the numerical order of accuracy using formula (\ref{orden}). The results for Gauss-Legendre formulas of two, three, four and five points are presented in Figure \ref{err_GL_composed} and Table \ref{tabla_exp4}. The dashed lines in Figure \ref{err_GL_composed} represent the theoretical decreasing of the error, that is $O(h^4)$ for the corrected two points formula, $O(h^6)$ for the corrected three points formula, $O(h^8)$ for the corrected four points formula, and $O(h^{10})$ for the corrected five points formula. These orders of accuracy can also be observed  in Table \ref{tabla_exp4}.

\begin{table}[!ht]
\begin{center}
\resizebox{18cm}{!} {
\begin{tabular}{|c|c|c|c|c|c|c|c|c|c|c|c|c|c|c|c|c|c|c|c|c|c|c|c|c|c|c|}
\hline $n=2^i$ & $2^3$& $2^4$  & $2^5$      & $2^6$     & $2^7$    & $2^8$ & $2^9$  %& $2^{10}$ %&$2^{11}$ &$2^{12}$&$2^{13}$%&$2^{14}$&$2^{15}$
      \\       
\hline
Error Gauss-Legendre 2 points ($E_i$)&5.18030e-02&3.60415e-02&8.54490e-03&8.49337e-03&2.12736e-03&2.13552e-03&5.33718e-04
      \\       
\hline 
$O_i$ &-&-0.52337&2.0765&0.0087262&1.9973&-0.0055253&2.0004
    \\
\hline     
Error Corr. Gauss-Legendre 2 points ($E_i$)& 1.94922e-02&1.33203e-03&8.52532e-05&5.36079e-06&3.35551e-07&2.09795e-08&1.31135e-09
      \\       
\hline     
$O_i$& -& 3.8712&3.9657&3.9912&3.9978&3.9995&3.9999
\\
\hline    
Error Gauss-Legendre 3 points ($E_i$)&4.14413e-02&1.33858e-02&1.04296e-02&3.31103e-03&2.60728e-03&8.30086e-04&6.51853e-04
      \\       
\hline 
$O_i$&- & 1.6304&0.36002&1.6553&0.34473&1.6512&0.34872
    \\
\hline     
Error Corr. Gauss-Legendre 3 points ($E_i$)& 1.63675e-04&2.87246e-06&4.62921e-08&7.29038e-10&1.14149e-11&1.81188e-13&7.10543e-15
      \\       
\hline     
$O_i$& -& 5.8324&5.9554&5.9886&5.997&5.9773&-
\\
\hline    
Error Gauss-Legendre 4 points ($E_i$)& 3.40572e-02&4.37203e-03&8.53484e-03&1.11711e-03&2.13317e-03&2.77810e-04&5.33339e-04
      \\       
\hline 
$O_i$&- & 2.9616&-0.96506&2.9336&-0.93323&2.9408&-0.94095
    \\
\hline     
Error Corr. Gauss-Legendre 4 points ($E_i$)& 8.88696e-07&3.99742e-09&1.62110e-11&6.75016e-14&3.55271e-15&3.55271e-15&1.06581e-14
      \\       
\hline     
$O_i$& -& 7.7965&7.9459&7.9078&4.2479&-&-
\\
\hline    
Error Gauss-Legendre 5 points ($E_i$)&1.46551e-02&1.39667e-02&3.58645e-03&3.47561e-03&9.01746e-04&8.69921e-04&2.25120e-04
      \\       
\hline 
$O_i$&- & 0.069416&1.9614&0.045291&1.9465&0.051836&1.9502
    \\
\hline     
Error Corr. Gauss-Legendre 5 points ($E_i$)& 3.50054e-09&4.02878e-12&3.55271e-15&3.55271e-15&3.55271e-15&3.55271e-15&7.10543e-15
      \\       
\hline     
$O_i$& -& 9.763&10.147&-&-&-&-
\\
\hline    

\end{tabular}
}
\caption{Grid refinement analysis for the integral of the function shown in (\ref{func_comp}) using some composed Gauss-Legendre quadrature rules.
}\label{tabla_exp4}
\end{center}
\end{table}

\section{Conclusions}\label{conc}

This study has introduced a new method to improve the accuracy of standard quadrature rules by incorporating additional correction terms. This strategy is especially effective when we have knowledge of a function's discontinuity location and the jumps in the function and its derivatives at that point.
The proposed correction terms, which are straightforward and easy to compute, have shown a considerable improvement over classical techniques. They allow the integration rule to be precise for piece-wise polynomials of degree $l$ or lower, thereby enhancing the accuracy of numerical integration.
This approach is particularly interesting for Gauss-Legendre quadrature formulas that employ a large number of points. In case of a discontinuity, the common adjustment would be to divide the integration interval using as a division point the position of the discontinuity, effectively doubling the number of points used. However, this can substantially escalate the computational cost if the function contains several isolated discontinuities. By utilising the proposed correction terms, we can retain the accuracy of the integration without the necessity to augment the number of points.
Thus, the proposed technique improves the accuracy of standard quadrature rules which enhances the accuracy when discontinuities are present, offering an efficient and easy-to-implement solution for numerical integration of piece-wise defined functions. %This advancement paves the way for more accurate and efficient computations in various scientific and engineering fields. Future endeavors will delve into further applications and potential enhancements of this method.

\section*{Acknowledgements}
Dr. Juan Ruiz has been supported by the Spanish national research project PID2019-108336GB-I00. The author Samala Rathan is supported by  NBHM, DAE, India (Ref. No. 02011/46/2021 NBHM(R.P.)/R \& D II/14874) and  IRG, IIPE Visakhapatnam with Grant No. IIPE/DORD/IRG/001. The author Shipra Mahata is supported by  UGC, India (NTA Ref No.- 211610046691). Dr. Dionisio F. Y\'a\~nez has been supported by project CIAICO/2021/227 (funded by Conselleria de Innovaci\'on, Universidades, Ciencia y Sociedad digital, Generalitat Valenciana) and by grant PID2020-117211GB-I00 (funded by MCIN/AEI/10.13039/501100011033).

\end{document}